\documentclass[10pt,leqno]{article}
\usepackage{amsfonts}
\usepackage{listings}
\usepackage{fullpage}
\usepackage{amssymb}
\usepackage{amsmath}
\usepackage{amsthm}
\newtheorem{theorem}{Theorem}

\newtheorem{corollary}{Corollary}

\newtheorem{lemma}{Lemma}

\begin{document}
\centerline{\Large \bf Density of integral sets with missing differences}\vskip 5mm

 \centerline{\large Quan-Hui Yang$^1$ and Min Tang$^2$}
 \footnotetext[1]{Supported by the Project of Graduate Education Innovation of Jiangsu Province
(CXZZ12-0381). Email:
yangquanhui01@163.com} \footnotetext[2]{Supported by the National Natural Science Foundation of China, Grant
No.10901002 and Anhui Provincial Natural Science Foundation, Grant No.1208085QA02.  Email:
tmzzz2000@163.com}
\begin{center}
{ \small 1. School of Mathematical Sciences, Nanjing Normal University, Nanjing
210023, China\\
2. Department of Mathematics, Anhui Normal University, Wuhu
241003, China} \end{center}

\vskip 3mm

\begin{abstract} Motzkin posed the problem of finding the maximal density $\mu(M)$
of sets of integers in which the differences given by a set $M$ do not occur. The problem
is already settled when $|M|\leq 2$ or $M$ is a finite arithmetic progression. In this paper,
we determine $\mu(M)$ when $M$ has some other structure. For example, we determine $\mu(M)$
when $M$ is a finite geometric progression.




\end{abstract}

\section{Introduction}

Let $\mathbb{N}$ be the set of all nonnegative integers. For a positive real number $x$ and
$S\subseteq \mathbb{N}$, we denote by $S(x)$ the number of elements $n\in S$ such that $n\leq x.$
The upper and lower densities of $S$, denoted by $\overline{\delta}(S)$ and $\underline{\delta}(S)$ respectively, are given by
$$\overline{\delta}(S):=\limsup_{x\rightarrow \infty}\frac{S(x)}{x}, \quad
\underline{\delta}(S):=\liminf_{x\rightarrow \infty}\frac{S(x)}{x}.$$
If $\overline{\delta}(S)=\underline{\delta}(S)$, we denote the common value by $\delta(S)$, and say
that $S$ has density $\delta(S).$

Given a set of positive integers $M$, we call a set $S\subseteq \mathbb{N}$ is an {\em $M$-set} if $a\in S,b\in S$
implies $a-b\notin M.$ In an unpublished problem collection, Motzkin \cite{Motzkin} posed the problem
of determining the quantity $$\mu(M): = \sup_S \overline{\delta}(S),$$ where the supremum is taken over all
$M$-sets $S$. In \cite{JCTA73}, Cantor and Gordon proved that if $|M|=1,$ then $\mu(M)=1/2$ and that if $M=\{m_1,m_2\},$
then $\mu(M)=[(m_1+m_2)/2]/(m_1+m_2).$ The following result is also proved.

\noindent{\bf  Theorem A.} {\em  Let $M_1=\{m_1,m_2,\ldots\}$ and $M_2=\{dm_1,dm_2,\ldots \},$ where $d$ is a positive integer.
Then $\mu(M_1)=\mu(M_2).$}

By Theorem A, we may assume that $\gcd(m_1,m_2,\ldots)=1$ for the purpose of determining
$\mu(M)$. Later, Haralambis \cite{JCTA77} determined $\mu(M)$ for most members of the families $\{1,j,k\}$ and
$\{1,2,j,k\}$. In 1999, Gupta and Tripathi \cite{Acta99}
completely determined $\mu(M)$ when $M$ is a finite arithmetic progression.

\noindent{\bf  Theorem B.} {\em  If $M=\{a,a+d,a+2d,\ldots,a+(n-1)d\}$ with $\gcd(a,d)=1$ and $n>1$,
then $$\mu(M)=\begin{cases} \frac{2a+(n-1)(d-1)}{2\{2a+(n-1)d\}} ~~~~~ \text{if} ~~ d \text{~is~odd}; \\
\frac 12   ~~~~~~~~~~~~~~~~~~\text{if} ~~ d \text{~is~even}.\\
\end{cases}$$}

In 2011, Pandey and Tripathi \cite{JNT11}
investigated this quantity when $M$ is related to arithmetic progressions.
For related results, one may refer to \cite{JCTA00},~\cite{Pandey} and \cite{JIS11}.

Motzkin's problem has connections with some other problems, such as the {\em $T$-colouring problem},
problems related to {\em the fractional chromatic number of distance graphs} and the {\em Lonely Runner Conjecture}. One may refer to
\cite{Bienia}, \cite{chang2}, \cite{Wu}.

Recently, Chen and Yang \cite{chen}, Khovanova and Konyagin \cite{Konyagin} studied the upper density
among sets of nonnegative integers in which no two elements have quotient belonging to $M$.

In this paper, the following results are proved.
\begin{theorem}\label{mainthm} Let $M=\{1,q,m_3,m_4,\ldots,m_s\}$, where $1<q<m_3<m_4<\cdots<m_s$.
If $m_i\equiv \pm 1 \pmod {q+1}$ for all integers $i\in \{3,4,\ldots,s\}$, then we have
$$\mu(M)=\begin{cases} 1/2 ~~~~~~ \text{if} ~~ q \text{~is~odd}; \\
\frac{q}{2(q+1)}   ~~~\text{if} ~~ q \text{~is~even}.\\
\end{cases}$$
\end{theorem}
From Theorem A and Theorem \ref{mainthm} we obtain the following corollary.

\begin{corollary}\label{coro1} If $M=\{a,aq,\ldots,aq^n\},$ where $a,q,n$ are positive integers with $q\geq 2$,
then $\mu(M)=1/2$ if $q$ is odd, and $\mu(M)=\frac{q}{2(q+1)}$ if $q$ is even.
\end{corollary}

In the next theorems, we shall consider some other sets $M$ with special structure.
\begin{theorem}\label{thm2} Suppose that $M=\{m_1,m_2,\ldots,m_n\}$ with $m_1<m_2<\cdots<m_n$ and satisfy the following two conditions:

 (i) $\{m_j-m_i:2\leq i<j\leq n\}\subseteq M;$

 (ii) the set $M$ does not contain a multiple of $n$.

\noindent { Then we have $\mu(M)= 1/n.$}
\end{theorem}

\begin{theorem}\label{thm3} Let $M=\{i+kj:0\leq k\leq n-1\}\cup \{j\}$ with $\gcd(i,j)=1$ and $n\geq 1$. Let $i+nj\equiv r\pmod {n+1}$ with $0\leq r\leq n$. If
$\gcd(r,i+nj)=1$, then $\mu(M)\geq \frac{i+nj-r}{(n+1)(i+nj)}$. Furthermore, if $r=n$, then $\mu(M)\geq \frac{i+nj-n}{(n+1)(i+nj)}$.
\end{theorem}

By Theorem \ref{thm2} and Theorem \ref{thm3}, we have the following corollary.

\begin{corollary}\label{cor2} If $M=\{i,j,i+j,i+2j\}$ with $i< j$ and $\gcd(i,j)=1,$ then
$$\mu(M)\begin{cases} =1/4 ~~~~~~ \text{if} ~~ i+3j\equiv 0~\text{or}~2~ \pmod 4; \\
\geq \frac{i+3j-r}{4(i+3j)}   ~~~\text{if} ~~ i+3j\equiv
r~\pmod 4,~\text{where}~r=\text{1 or 3}.
\end{cases}$$

\end{corollary}

\begin{theorem}\label{thm4} If $M=\{i+k_1j:0\leq k_1\leq n-1\}\cup \{k_2j:1\leq k_2\leq n\}$ with $n\geq 1$,
$\gcd(i,j)=1$ and $i+nj\equiv 1 \pmod {n+1},$ then
$$\mu(M)=\frac{i+nj-1}{(n+1)(i+nj)}.$$
\end{theorem}

\section{Preliminary Lemmas}\label{proofs}
In this section, we state two useful lemmas which give the lower and upper bound for $\mu(M)$.
\begin{lemma} (See \cite[Theorem 1]{JCTA73}.)\label{1} Let $M=\{m_1,m_2,m_3,\ldots\},$ and let $c$ and $m$ be positive integers such that $\gcd(c,m)=1.$
Put $$d=\min_{k}|cm_k|_m,$$
where $|x|_m$ denotes the absolute value of the absolutely least remainder of $x\pmod m$. Then $\mu(M)\geq d/m.$
\end{lemma}

\begin{lemma} (See \cite[Lemma 1]{JCTA77}.)\label{3} Let $M$ be a given set of positive integers, $\alpha$ a real number in the
interval $[0,1],$ and suppose that for any $M$-set $S$ with $0\in S$ there exists a positive integer $k$ (possibly dependent on
$S$) such that $S(k)\leq (k+1)\alpha$. Then $\mu(M)\leq \alpha.$
\end{lemma}

\section{Proof of Theorem \ref{mainthm}}

Suppose that $q$ is odd. For any $M$-set $S$, by $1\in M,$ we know that it contains no two consecutive integers. Thus, $\mu(M)\leq 1/2.$
On the other hand, by $m_i\equiv \pm 1 \pmod {q+1}$ and $2~|~q+1,$ we know that $M$ consists of only odd numbers.
Hence, the example $S=\{0,2,4,\ldots\}$ shows that equality can hold, and so $\mu(M)=1/2.$

Now we consider the case in which $q$ is even. If $\mu(M)>\frac{q}{2(q+1)}$, then
there exists an $M$-set $S$ and an interval $[c,c+q]$ such that $|S\cap [c,c+q]|>q/2$.
That is, $|S\cap [c,c+q]|\geq q/2+1.$ Noting that $1\in M$, we know that $S$ does not contain consecutive integers,
and so $S\cap [c,c+q]=\{c,c+2,\ldots,c+q\}$. It follows that $q\in S-S$, a contradiction. Hence, $\mu(M)\leq  \frac{q}{2(q+1)}$.

Next we shall prove that $\mu(M)\geq  \frac{q}{2(q+1)}$. Since $$\frac{q}{2}\times m_i \equiv \pm \frac{q}{2} \pmod {q+1}$$
for all integers $i\in \{3,4,\ldots,s\}$ and $1\cdot \frac q2 \equiv \frac q2 \pmod {q+1},$ $q\cdot \frac q2\equiv -\frac q2 \pmod {q+1},$ by Lemma \ref{1}, taking $c=q/2$ and $m=q+1,$ we have $\gcd(c,m)=1$ and $d=q/2.$ Thus,
$\mu(M)\geq \frac{q}{2(q+1)}.$

Therefore, $\mu(M)=\frac{q}{2(q+1)}$ if $q$ is even.

\section{Proof of Theorem \ref{thm2}}

For any positive integer $x$ and an $M$-set $S\subseteq [0,x]$,
we shall prove that $|S|\leq (x+m_n+1)/n.$

First, we prove $|S+M|\geq (n-1)|S|$ by induction on $|S|$. Clearly, it is true for $|S|=1$. Now suppose that $|S'+M|\geq (n-1)|S'|$ holds for all $M$-sets $S'$ with $|S'|<|S|.$ Write $S=\{b_1,b_2,\ldots,b_{|S|}\}.$ By $\{m_j-m_i:2\leq i<j\leq n\}\subseteq M,$ it follows that
$$\left( b_{|S|}+\{m_2,m_3,\ldots,m_n\} \right) \cap  \left(\{b_1,b_2,\ldots,b_{|S|-1}\}+M \right)=\emptyset.$$
Otherwise, there exist three integers $i,j,k$ with $2\leq i< j\leq n$ and $1\leq k \leq |S|-1$ such that
$b_{|S|}+m_i=b_k+m_j$, and then $b_{|S|}-b_k=m_j-m_i\in M$, a contradiction.  Hence, by the induction hypothesis, we have
$$|\{b_1,\ldots,b_{|S|}\}+M|\geq (n-1)+|\{b_1,\ldots,b_{|S|-1}\}+M|\geq (n-1)|S|.$$
By $S\cap (S+M)=\emptyset$ and $S\cup (S+M)\subseteq [0,x+m_n]$, it follows that
$n|S|\leq |S|+|S+M|\leq x+m_n+1$, and so $|S|\leq (x+m_n+1)/n.$

Hence, $\mu(M)=\sup_S \overline{\delta}(S)\leq \lim_{x\rightarrow \infty}(x+m_n+1)/(nx)=1/n.$

On the other hand, since $M$ does not contain a multiple of $n$, the set $\{0,n,2n,\ldots\}$ is an $M$-set. So $\mu(M)\geq 1/n$.

Therefore, $\mu(M)=1/n.$

\section{Proof of Theorem \ref{thm3}}

Let $t=i+nj$. Then $t\equiv r \pmod {n+1}.$ We consider the following two cases.

{\bf Case 1:}
$\gcd(r,t)=1.$ By $\gcd(i,j)=1,$ we have $\gcd(j,t)=1.$
Then there exists an integer $x$ such that $$xj\equiv \frac{t-r}{n+1} \pmod t.$$ Since
$\gcd(r,t)=1$, it follows that $\gcd(\frac{t-r}{n+1},t)=1$, and then $\gcd(x,t)=1.$ For such $x$ we have
$$x(i+kj)\equiv \frac{(k+1)(t-r)}{n+1}+r \pmod t$$
for $k=0,1,\ldots,n-1.$
Noting that $\gcd(x,t)=1$ and $$x(i+(n-1)j)\equiv \frac{n(t-r)}{n+1}+r\equiv -\frac{(t-r)}{n+1}\pmod t,$$ by Lemma \ref{1}, we have $\mu(M)\geq \frac{t-r}{(n+1)t}=\frac{i+nj-r}{(n+1)(i+nj)}.$

{\bf Case 2:} $r=n$ and $\gcd(r,t)>1.$ Then there exists an integer $x'$ such that $$x'j\equiv \frac{t+1}{n+1} \pmod t.$$
Since $\gcd(\frac{t+1}{n+1},t)=1$, we have $\gcd(x',t)=1.$
For such $x'$ we have
$$x'(i+kj)\equiv \frac{(k+1)(t-n)}{n+1}+k \pmod t$$
for $k=0,1,\ldots,n-1.$
Noting that $\gcd(x',t)=1$,
$$x'i\equiv \frac{t-n}{n+1} \pmod t,$$ and $$x'(i+(n-1)j)\equiv \frac{n(t-n)}{n+1}+n-1\equiv -\frac{(t+1)}{n+1}\pmod t,$$ by Lemma \ref{1}, we have $\mu(M)\geq \frac{t-n}{(n+1)t}=\frac{i+nj-n}{(n+1)(i+nj)}.$

\section{Proof of Corollary \ref{cor2}}
We consider the following three cases.

{\bf Case 1:} $i+3j\equiv r \pmod 4$ with $r=1$ or $3$. By Theorem \ref{thm3}, we have $\mu(M)\geq \frac{i+3j-r}{4(i+3j)}$.

{\bf Case 2:} $i+3j\equiv 0\pmod 4$, then we have $i\equiv j\equiv
1\pmod 4$ or $i\equiv j\equiv 3\pmod 4$. In this case, the set $M$ does
not contain a multiple of $4$, thus by Theorem \ref{thm2}, we have
$\mu(M)=1/4$.

{\bf Case 3:} $i+3j\equiv 2\pmod 4$, then we have $i\equiv 1\pmod
4$, $j\equiv 3\pmod 4$ or $i\equiv 3\pmod 4$, $j\equiv 1\pmod 4$.
In this case, we know that the set $M$ contains only one even
number $i+j$. Let
$$S=\bigcup_{k=0}^{\infty}\big(2k(i+j)+\{0,2,4,\ldots,i+j-2\}\big).$$
Clearly, $S$ is an $M$-set and $\delta(S)=1/4$. Hence $\mu(M)\geq
1/4.$ By the proof of Theorem \ref{thm2}, we have $\mu(M)\leq 1/4.$ Therefore,
$\mu(M)=1/4.$

\section{Proof of Theorem \ref{thm4}}

First we follow the proof of Theorem \ref{thm3} and show that $\mu(M)\geq \frac{i+nj-1}{(n+1)(i+nj)}$. Let $i+nj=t$. Then
there exists an integer $x$ such that
$$xj\equiv \frac{t-1}{n+1} \pmod t.$$ Clearly, $\gcd(x,t)=1.$
For such $x$ we have
$$x(i+k_1j)\equiv \frac{(k_1+1)(t-1)}{n+1}+1 \pmod t$$ for
$k_1=0,1,\ldots,n-1$ and $$x(k_2j)\equiv \frac{k_2(t-1)}{n+1} \pmod t$$ for
$k_2=1,2,\ldots,n.$
Noting that
$$x(i+(n-1)j)\equiv \frac{n(t-1)}{n+1}+1\equiv -\frac{t-1}{n+1}  \pmod t$$
and $$x(nj)\equiv \frac{n(t-1)}{n+1}\equiv -\frac{t-1}{n+1}-1 \pmod t,$$
by Lemma \ref{1} we have $\mu(M)\geq \frac{t-1}{(n+1)t}$.

Now we will prove $\mu(M)\leq \frac{t-1}{(n+1)t}$. Let $S$ be any $M$-set with $0\in S$. Then
for $t=i+nj$, $$\bigcup_{m=1}^{(t-n-2)/(n+1)}A_m \cup B$$ is a decomposition of $\{0,1,2,\ldots,t-1\}$ into
disjoint sets, where $$B=\{0,j\}\cup \{i+kj:0\leq k\leq n-1\},$$ $$A_m=\{(j-i)m+i+kj:1\leq k\leq n-1\}\cup\{(j-i)m,(j-i)m+j\}$$
with $1\leq m\leq (t-n-2)/(n+1)$
and the elements of each set $A_m$ are considered modulo $t$.

Since $S$ is an $M$-set and $0\in S$, it follows that $|A_m\cap S|\leq 1$ for each $m$ and $|S\cap B|=1.$

Hence, $S(t-1)\leq 1+(t-n-2)/(n+1)=(t-1)/(n+1)$ for any $M$-set $S$.
By Lemma \ref{3}, it follows that $\mu(M)\leq (t-1)/((n+1)t).$

Therefore, we obtain $$\mu(M)= \frac{t-1}{(n+1)t}=\frac{i+nj-1}{(n+1)(i+nj)}.$$

\section{Acknowledgments}  
We sincerely thank the anonymous
referee for his/her detailed comments. 

\end{document}